\documentstyle[12pt]{article}

\large\normalsize
\newcommand{\eqref}[1]{(\ref{#1})}
\newtheorem{proposition}{Proposition}[section]

\newtheorem{remark and definition}[proposition]{Remark and definition}
\newtheorem{theorem}[proposition]{Theorem}
\newtheorem{important example}[proposition]{Exemple important}
\newtheorem{principal mutations}[proposition]{Mutations principales}
\newtheorem{corollary}[proposition]{Corollary}
\newtheorem{remark}[proposition]{Remark}
\newtheorem{particular case}[proposition]{Cas particulier}

\newtheorem{notation and remark}[proposition]{Notation and remark}
\newtheorem{remarks}[proposition]{Remarks}
\newtheorem{note and remark}[proposition]{Note and remark}
\newtheorem{lemma}[proposition]{Lemma}

\newtheorem{examples}[proposition]{Examples}
\newtheorem{definitions}[proposition]{Definitions}
\newtheorem{definitions and notations}[proposition]{Definitions and notations}
\newtheorem{principal isotopes}[proposition]{Principal isotopes}

\newtheorem{problem}[proposition]{Problem}

\newtheorem{third power-associative algebra}[proposition]{Third power-associative algebra}
\newtheorem{remarks and open problem}[proposition]{Remarques et probl\`eme ouvert}

\def\ait{{ \mathbb{A}}}
\def\ait{\hbox{\it l\hskip -2pt A}}
\def\cit{\hbox{\it l\hskip -5.5pt C\/}}
\def\oit{\hbox{\it l\hskip -5.5pt O\/}}
\def\rit{{ \mathbb{R}}}
\def\rit{\hbox{\it l\hskip -2pt R}}

\def\pit{\hbox{\it l\hskip -2pt P}}
\def\nit{\hbox{\it l\hskip -2pt N}}

\def\hit{\hbox{\it l\hskip -2pt H}}

\begin{document}

\sloppy
\renewcommand{\thesection}{\arabic{section}}
\renewcommand{\theequation}{\thesection.\arabic{equation}}

\setcounter{page}{1}

\title{\bf \sf On finite-dimensional absolute-valued algebras satisfying $(x^p,x^q,x^r)=0$}
\author{A. Chandid, M. I. Ram\'irez and A. Rochdi}
\maketitle

\vspace{0.1cm}
\begin{abstract} By means of principal isotopes $\hit(a,b)$ of algebra $\hit$ \ [Ra 99]
we give an exhaustive description of all $4$-dimensional absolute
valued algebras satisfying $(x^p,x^q,x^r)=0$ for fixed integers
$p,q,r\in\{1,2\}.$ For such an algebras the number $N(p,q,r)$ of
isomorphisms classes is either finite included between two and
three, or is infinite. Concretely
\begin{enumerate} \item $N(1,1,1)=N(1,1,2)=N(1,2,1)=N(2,1,1)=2,$ \item $N(1,2,2)=N(2,2,1)=3,$
\item $N(2,1,2)=N(2,2,2)=\infty.$
\end{enumerate} Besides, each one of above algebras contains $2$-dimensional subalgebras. However,
the problem in dimension $8$ is far from being completely solved.
In fact, there are $8$-dimensional absolute-valued algebras,
containing no $4$-dimensional subalgebras, satisfying
$(x^2,x,x^2)=(x^2,x^2,x^2)=0.$
\end{abstract}

\vspace{0.1cm} {\bf Keywords.} Absolute valued algebra, central
(flexible) idempotent, principal isotopes of $\hit.$

\vspace{0.5cm}
\section{Introduction}

\vspace{0.2cm} \hspace{0.3cm} Absolute-valued algebras, introduced
in 1918 [Ost 18], constitute nowadays a distinguished and
attractive class of non-associative algebras, by its widely
observed variety, and we find in [Rod 04] a perfect compilation of
the results around the theory, previous in 2004. During the last
decade several works on finite-dimensional absolute-valued
algebras appeared [Ra 99], [Roc 03], [Rod 04], [CM 05], [CR 08],
[RR 09], [CKMMRR]. A fundamental result on this subject, given by
Albert [A 47], asserts that the dimension of every
finite-dimensional absolute-valued algebra $A$ is $1,2,4$ or $8,$
and $A$ is isotopic to one of classical (unital) absolute-valued
algebras \ $\rit,$ $\cit,$ $\hit$ (real quaternion algebra) or
$\oit$ (real octonion algebra). It follows easily that \ $\rit,$
$\cit,$ $^*\cit,$ $\cit^*,$ $\stackrel{*}{\cit}$ are the unique
absolute-valued real algebras of dimension $\leq 2.$ Using the
so-called principal isotopes of $\hit$ formed by four families
depending on a pair of norm-one elements of $\hit,$ Ram\'irez [Ra
99] classified all $4$-dimensional absolute-valued real algebras
and solved the isomorphism problem. Calder\'on and Mart\'in [CM
05] gave a refinement of above classification. Recently a
description of all $8$-dimensional absolute-valued algebras was
given [CKMMRR]. On the other hand, classical results specify those
absolute-valued algebras which satisfy to a familiar identity as
associativity [Ost 18], [M 38] or commutativity [UW 60] and where
extended to power-associativity [EM 80] and flexibility [EM 81]
respectively. Concretely

\begin{enumerate} \item Every absolute-valued associative real algebra is
isometrically isomorphic to either \ $\rit,$ $\cit$ or $\hit.$
\item $\rit,$ $\cit,$ and $\stackrel{*}{\cit}$ \ are the unique
absolute-valued commutative real algebras. \item $\rit,$ $\cit,$
$\hit,$ and $\oit$ \ are the unique absolute-valued
power-associative real algebras. \item Every absolute-valued
flexible real algebra is finite-dimensional.
\end{enumerate}

\vspace{0.1cm} \hspace{0.3cm} By considering third
power-associative algebras, that is algebras satisfying
$(x,x,x)=0$ where $(.,.,.)$ means associator, El-Mallah [Elm 87]
and, independently, Elduque-P\'erez [EP 94] gave a
finite-dimensional extension of above results by showing that
every finite-dimensional absolute-valued algebra satisfying
$(x,x,x)=0$ is flexible and equal to either \ $\rit,$ $\cit,$
$\stackrel{*}{\cit},$ $\hit,$ $\stackrel{*}{\hit},$ $\oit,$
$\stackrel{*}{\oit},$ or the algebra \ $\pit$ \ of
pseudo-octonions. In a general way absolute-valued algebras
satisfying to an identity of the form $(x^p,x^q,x^r)=0$ for fixed
integers $p,q,r\in\{1,2\},$ including $(x,x,x)=0,$ where
extensively studied [Elm 82, 83, 87, 88, 97, 01], [EP 94], [EE
04], [Cu 06], [CR 08], [RR 09], [MRR]. For such an algebra $A,$
the existence of either a non-zero central idempotent [Elm 83,
01], [EE 04], [Cu 06], [CR 08], [RR 09], or a left-unit [CR 08]
carry away finite dimensionality. Moreover, a complete description
of these algebras was given in both cases ([Elm 88], [Cu 06], [CR
08], [RR 09]) and ([CR 08], [MRR]), respectively. However there
are infinite-dimensional absolute-valued algebras satisfying to
the identity $(x^2,y,x^2)=0,$ in two variables [CR 08]. In the
other hand it seems that, with the exception of identity
$(x,x,x)=0,$ no previous study took place for an identity of the
form $(x^p,x^q,x^r)=0$ over an arbitrary absolute-valued algebra
of finite dimension $\geq 4.$ The aim of this paper is to make
some contributions to this subject. Here we give a description of
all four-dimensional absolute-valued algebras satisfying an
identity of the form $(x^p,x^q,x^r)=0$ for fixed integers
$p,q,r\in\{1,2\}$ (Theorem {\bf 3.24}). We show the following main
result:

\vspace{0.6cm} {\bf Theorem.} {\em Let $A$ be a four-dimensional
absolute valued algebra. Then $A$ satisfies to an identity of the
form $(x^p,x^q,x^r)=0,$ for fixed integers $p,q,r\in\{1,2\},$ if
and only if $A$ is isomorphic to one of principal isotopes of
$\hit$ described by the following table:}

\[ \begin{tabular}{cc} \\ \hline
\multicolumn{1}{|c|}{$A$ satisfies to} & \multicolumn{1}{|c|}{The list of isomorphism classes} \\
\hline \multicolumn{1}{|c|}{$(x,x,x)=0$ or $(x,x^2,x)=0$} &
\multicolumn{1}{|c|}{$\hit,$ $\stackrel{*}{\hit}$} \\
\hline \multicolumn{1}{|c|}{$(x,x,x^2)=0$} &
\multicolumn{1}{|c|}{$\hit,$ $\hit^*$} \\
\hline \multicolumn{1}{|c|}{$(x,x^2,x^2)=0$} &
\multicolumn{1}{|c|}{$\hit,$ $\hit^*,$ $\hit^*(1,i)$} \\
\hline \multicolumn{1}{|c|}{$(x^2,x,x)=0$} &
\multicolumn{1}{|c|}{$\hit,$ ${^*\hit}$}
\\ \hline \multicolumn{1}{|c|}{$(x^2,x^2,x)=0$} &
\multicolumn{1}{|c|}{$\hit,$ ${^*\hit},$ ${^*\hit}(i,1)$}
\\ \hline \multicolumn{1}{|c|}{$(x^2,x,x^2)=0$} &
\multicolumn{1}{|c|}{$\hit,$ $\stackrel{*}{\hit}$} \\
\multicolumn{1}{|c|}{} &
\multicolumn{1}{|c|}{${^*\hit}(1,e^{i\alpha}),$
$\hit^*(e^{i\alpha},1)$ with $\alpha\in [0,\pi[$} \\ \hline
\multicolumn{1}{|c|}{} & \multicolumn{1}{|c|}{$\hit,$
$\stackrel{*}{\hit}$}
\\ \multicolumn{1}{|c|}{$(x^2,x^2,x^2)=0$} &
\multicolumn{1}{|c|}{${^*\hit}(i,1),$ ${^*\hit}(i,i),$
$\hit^*(1,i),$
$\hit^*(i,i)$} \\
\multicolumn{1}{|c|}{} &
\multicolumn{1}{|c|}{${^*\hit}(1,e^{i\alpha}),$
$\hit^*(e^{i\alpha},1) \mbox{ with } \alpha\in [0,\pi[$}
\\ \hline
\end{tabular}
\]

\vspace{0.5cm} \hspace{0.3cm} We can easily see that each of above
algebras contains $2$-dimensional subalgebras. However, the
problem in dimension $8$ is very hard. There are $8$-dimensional
absolute-valued algebra satisfying $(x^2,y,x^2)=0$ containing no
$4$-dimensional subalgebras. It is well known that every third
power-associative algebra, over a field of characteristic $\neq 2$
having at least three elements, satisfies to $(x,x^2,x)=0.$ In the
other hand we see that, among those $4$-dimensional
absolute-valued algebras, the identity $(x,x^2,x)=0$ carry away
third power-associativity. However the question to know if an
$8$-dimensional absolute-valued algebra which satisfies both
identities $(x,x^2,x)=(x^2,y,x^2)=0$ is third power-associative,
seems still to be an open problem.

\vspace{1cm} \section{Notations and preliminary results}

\vspace{0.5cm} \begin{definitions and notations} An algebra over
an arbitrary field $K$ is a vector space over $K$ endowed with a
bilinear mapping \ $(x,y)\mapsto xy$ \ from $A\times A$ to $A$
called the product of the algebra $A.$ For $x,y,z$ in $A,$ we
denote by $(x,y,z)$ the associator \ $(xy)z-x(yz)$ \ of $x,y,z$
and by $[x,y]$ the commutator \ $xy-yx$ \ of $x,y.$ We denote by
$L_x,$ $R_x$ the linear operators of $A$ defined by \ $y\mapsto
xy,$ \ $y\mapsto yx$ \ respectively. An element $e\in A$ is said
to be {\bf central} if it satisfies $[e,x]=0$ for all $x\in A,$
or, in other word $L_e=R_e.$ The element $e\in A$ is said to be
{\bf flexible} if it satisfies $(e,x,e)=0$ for all $x\in A,$ or,
in other words $L_e\circ R_e=R_e\circ L_e.$
\begin{enumerate} \item The
algebra $A$ is said to be \begin{enumerate} \item {\bf Third
power-associative} if it satisfies the identity $(x,x,x)=0.$ \item
{\bf Alternative} if it satisfies both identities
$(x,x,y)=(y,x,x)=0.$ A well known Artin's theorem asserts that an
algebra is alternative if and only if the subalgebra generated by
two arbitrary elements is associative {\em ([Sc 66] Theorem {\bf
3.1} p. 29)}. \item {\bf Flexible} it satisfies the identity
$(x,y,x)=0.$ \item {\bf Quadratic} if it contains an unit element
$e$ and $e,x,x^2$ are linearly dependent for all $x\in A.$ The
elements in the set

\[ Im(A)=\{x\in
A:x^2\in Ke \mbox{ and } x\notin Ke-\{0\}\} \]

are called the purely imaginary elements in $A.$ If, moreover, the
characteristic of $K$ is not two, then $Im(A)$ is a subspace of
$A$ supplementary to $Ke,$ i.e., $A=Ke\oplus Im(A)$ {\em ([HKR 91]
p. 227-228)}.
\end{enumerate} \item A linear operation \ $*:x\mapsto x^*$ defined
on algebra $A$ is said to be an involution if it satisfies, for
all $x,y\in A,$ the following conditions
\begin{enumerate}
\item $(x^*)^*=x,$ \item $xx^*=x^*x,$ \item $(xy)^*=y^*x^*.\Box$
\end{enumerate}
\end{enumerate}
\end{definitions and notations}

\vspace{0.2cm}
\begin{proposition} Let $A$ be an algebra over a field $K$ of characteristic $\neq 2$
containing at least three elements. If algebra $A$ is third
power-associative then $A$ satisfies to the identity
$(x,x^2,x)=0.$
\end{proposition}

\vspace{0.1cm} {\bf Proof.} A simple linearization of $(x,x,x)=0,$
taking into account that $K$ contains two distinct non-zero
elements, gives \ $[x^2,y]+[xy+yx,x]=0.$ The result follows by
putting $y=x^2$ in above identity, taking into account the
characteristic of $K.\Box$

\vspace{0.2cm}
\begin{definitions} An algebra $A,$ over real numbers $\rit$ or
complex numbers $\cit,$ is said to be absolute valued if the space
$A$ is endowed with a norm \/ $||.||$ \/ satisfying \/
$||xy||=||x||\hspace{0.1cm}||y||$ \/ for all $x,y\in A.\Box$
\end{definitions}

\vspace{0.2cm} \begin{remarks} The algebras \ $\cit,$ $\hit$
(quaternions) and \ $\oit$ (octonions) can be built from \ $\rit$
\ by iterating the Cayley-Dickson doubling process {\em ([HKR 91]
p. 256-257)}. If $A$ stands for either $\rit,$ $\cit,$ $\hit$ or
$\oit,$ with unit element $1$ and standard involution (SI)
$\sigma_A:x\mapsto\overline{x},$ then
\begin{enumerate} \item $A$ is a
quadratic and alternative algebra, and the expression $\rit\oplus
Im(A)$ is nothing other than the decomposition of $A$ into a
direct sum of the eigenspaces $E_1(\sigma_A),$ $E_{-1}(\sigma_A)$
associated, respectively, to eignevalues $1$ and $-1.$ Thus every
element $x\in A$ can be written as a sum of its real part \
$\frac{1}{2}(x+\bar{x}):=Re(x)\in\rit$ and purely imaginary part \
$\frac{1}{2}(x-\bar{x}):=Im(x)\in Im(A).$ \item The usual inner
product \[ (.|.):A\times A\rightarrow\rit^+\hspace{0.3cm}
(x,y)\mapsto(x|y)=Re(x\bar{y}) \] {\em ([HKR 91] p. 208, 253)}
satisfies to \ $(xy|xy)=(x|x)(y|y)$ for all $x,y\in A$ and the
(euclidian) norm $||.||$ of $A,$ given by
$||x||=\sqrt{x\overline{x}}=\sqrt{\overline{x}x}$ for all $x\in
A,$ is an absolute-value. Thus $(A,||.||)$ is a pre-Hilbert real
absolute-valued algebra and the SI $\sigma_A$ of $A=\rit\oplus
Im(A)$ \ is a linear isometry. Besides, the following equalities
holds for all $x,y\in A:$
\begin{enumerate} \item $(x|y)=(\overline{x}|\overline{y}).$ This
follows from equality $(\overline{x}|\overline{x})=(x|x)$ by
linearization. \item $(.|.)$ is a trace form on $A$ {\em ([BBO 82]
Lemma {\bf 1.2})}:
\begin{eqnarray} (xy|z) &=& (x|yz)
\end{eqnarray} \item The triple
product identity {\em ([HKR 91] p. 209, 254)}:
\begin{eqnarray} yxy &=&
2(\overline{x}|y)y-||y||^2\overline{x}.
\end{eqnarray}
\end{enumerate}
\end{enumerate}
\end{remarks}

\vspace{5cm} \hspace{0.3cm} According to [A 47], two
absolute-valued algebras $A$ and $B,$ are said to be isotopic if
there exist linear isometries $\Phi_1, \Phi_2, \Phi_3$ from $A$
onto $B$ satisfying $\Phi_1(xy)=\Phi_2(x)\Phi_3(y)$ for all
$x,y\in A.$ One of fundamental results about finite-dimensional
absolute-valued algebras, proved by Albert ([A 47] Theorem {\bf
2}), is the following:

\vspace{0.2cm} \begin{theorem} Let $A$ be a finite-dimensional
absolute-valued real algebra. Then $A$ is isotopic to either \
$\rit,$ $\cit,$ $\hit$ or $\oit.$ Therefore $A$ has dimension
$1,2,4$ or $8$ and the absolute-value of $A$ comes from an inner
product.$\Box$
\end{theorem}

\vspace{0.2cm} \begin{examples} .
\begin{enumerate} \item {\em\bf (Standard isotopes of $\ait$).}
For $\ait$ equal to either $\cit,$ $\hit$ or $\oit$ with SI
$x\mapsto\overline{x},$ let us denote by \ $^*\ait,$ $\ait^*$ and
$\stackrel{*}{\ait}$ the absolute-valued real algebras obtained by
endowing the normed space $\ait$ with the products $x\odot
y=\overline{x}y,$ $x\odot y=x\overline{y}$ and $x\odot
y=\overline{x}\hspace{0.1cm}\overline{y},$ respectively. Algebras
$\ait,$ $^*\ait,$ $\ait^*,$ $\stackrel{*}{\ait}$ are called
standard isotopes of $\ait.$ It follows easily from Theorem {\bf
2.5} that $\cit,$ $^*\cit,$ $\cit^*$ and $\stackrel{*}{\cit}$ \
are the unique $2$-dimensional absolute-valued real algebras {\em
(see [Rod 94], [EP 99])}. \item {\em\bf (Principal isotopes of
$\hit$).} These are the absolute-valued algebras $\hit_1(a,b),$
$\hit_2(a,b),$ $\hit_3(a,b)$ and $\hit_4(a,b)$ obtained from fixed
norm-one elements $a,b$ in $\hit$ by endowing the normed space
$\hit$ with the products \begin{eqnarray*} x\stackrel{a,b}{\odot_1}x &=& axyb, \\
x\stackrel{a,b}{\odot_2}x &=& \overline{x}ayb, \\
x\stackrel{a,b}{\odot_3}y &=& axb\overline{y}, \\
x\stackrel{a,b}{\odot_4}y &=&
a\overline{x}\hspace{0.1cm}\overline{y}b. \end{eqnarray*}

These four types of absolute-valued algebras are called principal
isotopes of $\hit$ denoted $\hit(a,b),$ $^*\hit(a,b),$
$\hit^*(a,b)$ and $\stackrel{*}{\hit}(a,b),$ respectively {\em [Ra
99]}.$\Box$
\end{enumerate}
\end{examples}

\vspace{0.2cm} \hspace{0.3cm} M. I. Ram\'irez [Ra 99] gave a
classification of all $4$-dimensional absolute-valued algebras by
using Theorem {\bf 2.5} and the description of all linear
isometries on $\hit$ ([HKR 91] p. 215):

\vspace{0.2cm}
\begin{theorem} Every four-dimensional absolute valued real algebra is isomorphic to a principal
isotope of $\hit.$ Moreover two principal isotopes $\hit_n(a,b)$
and $\hit_{m}(a',b')$ of $\hit$ are isomorphic if and only if
$n=m$ and the equalities $a'p=\varepsilon pa$ and $b'p=\delta pb$
hold for some norm-one element $p\in\hit$ and some $\varepsilon,
\delta\in\{1,-1\}.\Box$
\end{theorem}

\vspace{0.2cm}
\begin{remarks} Let $n$ be in $\{1,2,3,4\}.$ \begin{enumerate} \item Algebras $\hit_n(-a,-b),$ $\hit_n(-a,b),$
$\hit_n(a,-b),$ $\hit_n(a,b)$ are isomorphic. \item A refinement
of Theorem {\bf 2.7} can be found in {\em [CM 05]}. \item We check
easily that, among those principal isotopes of $\hit,$ algebras
$\hit(a,\overline{a}),$ $\stackrel{*}{\hit}(a,\overline{a})$ are
the unique ones which contain a non-zero central idempotent $e,$
inevitably equal to $1.\Box$
\end{enumerate}
\end{remarks}

\vspace{0.2cm} \hspace{0.3cm} Ram\'irez [Ra 99] gave, among
others, a precise description of all $4$-dimensional
absolute-valued algebras which contain $2$-dimensional
sub-algebras:

\vspace{0.2cm}
\begin{theorem} A $4$-dimensional absolute-valued real algebra contains
subalgebras of dimension $2$ if and only if it is isomorphic to a
principal isotope of $\hit$ whose parameters $a,b$ are precisely
subjeted to the conditions given by the following table

\[ \begin{tabular}{ccccc} \\ \hline
\multicolumn{1}{|c|}{$\hit(a,b)$} &
\multicolumn{1}{c|}{${^*\hit}(a,b)$} &
\multicolumn{1}{c|}{$\hit^*(a,b)$} &
\multicolumn{1}{c|}{$\stackrel{*}{\hit}(a,b)$} &
\multicolumn{1}{|c|}{contains} \\
\hline \multicolumn{1}{|c|}{$ab=ba$} & \multicolumn{1}{c|}{$a,
ba\in Im(\hit)$} & \multicolumn{1}{c|}{$b, ba\in Im(\hit)$} &
\multicolumn{1}{c|}{$a, b\in Im(\hit)$} &
\multicolumn{1}{|c|}{$\cit$} \\
\multicolumn{1}{|c|}{$a, ab\in Im(\hit)$} &
\multicolumn{1}{c|}{$ab=ba$} & \multicolumn{1}{c|}{$a, ba\in
Im(\hit)$} & \multicolumn{1}{c|}{$a, ba\in Im(\hit)$} &
\multicolumn{1}{|c|}{${^*\cit}$} \\
\multicolumn{1}{|c|}{$b, ba\in Im(\hit)$} &
\multicolumn{1}{|c|}{$b, ba\in Im(\hit)$} &
\multicolumn{1}{c|}{$ab=ba$} & \multicolumn{1}{c|}{$b, ba\in
Im(\hit)$} &
\multicolumn{1}{|c|}{$\cit^*$} \\
\multicolumn{1}{|c|}{$a, b\in Im(\hit)$} &
\multicolumn{1}{|c|}{$a, b\in Im(\hit)$} &
\multicolumn{1}{|c|}{$a, b\in Im(\hit)$} &
\multicolumn{1}{c|}{$ab=ba$} & \multicolumn{1}{|c|}{$\stackrel{*}{\cit}$} \\
\hline
\end{tabular} \]
\end{theorem}

\vspace{0.3cm} \hspace{0.3cm} We shall need the following useful
result:

\vspace{0.2cm}
\begin{lemma} For $a,b$ in $\hit,$ with the same norm, the the following are equivalents:
\begin{enumerate}
\item $a$ and $b$ are conjugated. \item $Re(a)=Re(b).$
\end{enumerate}
\end{lemma}

\vspace{0.1cm} {\bf Proof.} The implication {\bf 1 $\Rightarrow$
2} is clear and it remains only to show {\bf 2 $\Rightarrow$ 1}
with $a,b\in Im(\hit).$ Indeed, we have:
$b(a+b)=ba+b^2=ba+a^2=(a+b)a.$
\begin{enumerate} \item If $b\neq -a,$ then $b=va\overline{v}$ where
$v=||a+b||^{-1}(a+b).$ \item If $b=-a\neq 0,$ then for all
norm-one $u\in Im(\hit),$ orthogonal to $a,$ we have:
$b=u^2a=-uau=ua\overline{u}.\Box$
\end{enumerate}

\vspace{1cm} \section{Classification in dimension four}

\vspace{0.8cm} \subsection{Principal isotopes of $\hit$ satisfying
to $(x^p,x^q,x^r)=0$}

\vspace{0.5cm} \hspace{0.3cm} Now we will study those
finite-dimensional absolute-valued algebras which satisfy an
identity of the form $(x^p,x^q,x^r)=0$ for fixed integers
$p,q,r\in\{1,2\}.$ The situation in dimension one, corresponding
to the associative algebra $\rit,$ where all the identities
$(x^p,x^q,x^r)=0$ are satisfied, is trivial. It is not so much in
dimension two:

\vspace{0.4cm} \begin{proposition} Among all $2$-dimensional
absolute-valued algebras
\begin{enumerate} \item $\cit$ is the unique which satisfies all the identities $(x^p,x^q,x^r)=0.$
\item The other algebras, $^*\cit,$ $\cit^*,$
$\stackrel{*}{\cit},$ satisfy to \ $(x^2,x,x^2)=(x^2,x^2,x^2)=0.$
Among these three last ones
\begin{enumerate} \item $^*\cit$ is the unique which satisfy to $(x^2,x,x)=0$ or else, to $(x^2,x^2,x)=0.$ \item
$\cit^*$ is the unique which satisfy to $(x,x,x^2)=0$ or else, to
$(x,x^2,x^2)=0.$ \item $\stackrel{*}{\cit}$ is the unique which
satisfy to $(x,x,x)=0$ or else, to $(x,x^2,x)=0.\Box$
\end{enumerate}
\end{enumerate}
\end{proposition}

\vspace{0.2cm} \hspace{0.3cm} We will show that the general
situation in dimension $4$ is very diversified and principal
isotopes of $\hit$ will be used for our study. First of all the
following remark, for standard isotopes of $\hit,$ will serve us
as guide:

\vspace{0.4cm} \begin{remark} Let us note that the associative
algebra $\hit$ satisfy to all the identities of the form
$(x^p,x^q,x^r)=0$ and that algebras $^*\hit,$ $\hit^*,$
$\stackrel{*}{\hit}$ satisfy to both identities
$(x^2,x,x^2)=(x^2,x^2,x^2)=0.$ That, among these last three
algebras
\begin{enumerate} \item $^*\hit$ is the unique
which satisfy to $(x^2,x,x)=0$ or else, to $(x^2,x^2,x)=0.$ \item
$\hit^*$ is the unique which satisfy to $(x,x,x^2)=0$ or else, to
$(x,x^2,x^2)=0.$ \item $\stackrel{*}{\hit}$ is the unique which
satisfy to $(x,x,x)=0$ or else, to $(x,x^2,x)=0.\Box$
\end{enumerate}
\end{remark}

\vspace{0.3cm} \hspace{0.3cm} We shall have the opportunity to use
the supplementary notations \/ $\hit_1=\hit,$ \/ $\hit_2=
{^*\hit},$ \/ $\hit_3=\hit^*,$ \/
$\hit_4=\hspace{0.1cm}\stackrel{*}{\hit}.$ In our computations the
first purely imaginary element $i$ in the canonical basis \
$\{1,i,j,k\}$ \ of algebra $\hit,$ will be used.

\vspace{0.5cm} \hspace{0.3cm} Now we need the following
preliminary result:

\vspace{0.3cm}
\begin{lemma} Assume that for $a\in\hit$ the equality \ $[axa,x]=0$ \ holds
for all $x\in Im(\hit).$ Then $a\in\rit.$
\end{lemma}

\vspace{0.1cm} {\bf Proof.} For all $x\in Im(\hit),$ we
have \begin{eqnarray*} 0 &=& [axa,x] \\
&=& [2(\overline{x}|a)a-||a||^2\overline{x},x] \hspace{0.2cm}
\mbox{ by triple
product identity } \\
&=& -2(x|a)[a,x]. \end{eqnarray*}

\vspace{0.2cm} If $a\notin\rit$ then the choice of an element
$x\in Im(\hit)$ neither orthogonal nor linearly dependent to
$Im(a)$ leads to a contradiction.$\Box$

\vspace{0.4cm} \hspace{0.3cm} Let us note now [pqrn] the equation
deduced from identity $(x^p,x^q,x^r)=0,$ if it occurs, in algebra
$\hit_n(a,b).$ Theses 32 equations are illustrated in the
following table:

\vspace{0.1cm}
\[ \begin{tabular}{ccc} \\ \hline
\multicolumn{1}{|c|}{Identity $(x^p,x^q,x^r)=0$} &
\multicolumn{1}{|c|}{Corresponding equation} &
\multicolumn{1}{|c|}{[pqrn]} \\
\hline \multicolumn{1}{|c|}{} &
\multicolumn{1}{|c|}{$ax^2bx=xax^2b$} &
\multicolumn{1}{|c|}{[1111]} \\
\multicolumn{1}{|c|}{$(x,x,x)=0$} & \multicolumn{1}{|c|}{$
\overline{b}\hspace{0.1cm}\overline{x}\hspace{0.1cm}\overline{a}xax=
\overline{x}a\overline{x}axb$} &
\multicolumn{1}{|c|}{[1112]} \\
\multicolumn{1}{|c|}{} &
\multicolumn{1}{|c|}{$axb\overline{x}b\overline{x}=xbx\overline{b}\hspace{0.1cm}\overline{x}\hspace{0.1cm}
\overline{a}$}
& \multicolumn{1}{|c|}{[1113]} \\
\multicolumn{1}{|c|}{} &
\multicolumn{1}{|c|}{$bx\overline{b}x^2=x^2\overline{a}xa$} &
\multicolumn{1}{|c|}{[1114]} \\ \hline \multicolumn{1}{|c|}{} &
\multicolumn{1}{|c|}{$ax^2bax^2=xaxax^2b$} &
\multicolumn{1}{|c|}{[1121]} \\
\multicolumn{1}{|c|}{$(x,x,x^2)=0$} &
\multicolumn{1}{|c|}{$\overline{b}\hspace{0.1cm}\overline{x}\hspace{0.1cm}\overline{a}xa\overline{x}ax=
\overline{x}a\overline{x}a\overline{x}axb$} &
\multicolumn{1}{|c|}{[1122]} \\
\multicolumn{1}{|c|}{} &
\multicolumn{1}{|c|}{$axb\overline{x}bx=xbaxb\overline{x}$} &
\multicolumn{1}{|c|}{[1123]} \\
\multicolumn{1}{|c|}{} &
\multicolumn{1}{|c|}{$\overline{b}x^2\overline{a}\hspace{0.1cm}\overline{b}x=\overline{x}\hspace{0.1cm}\overline{b}
a\overline{x}^2b$} & \multicolumn{1}{|c|}{[1124]} \\ \hline
\multicolumn{1}{|c|}{} &
\multicolumn{1}{|c|}{$axax^2b^2x=xa^2x^2bxb$} &
\multicolumn{1}{|c|}{[1211]} \\
\multicolumn{1}{|c|}{$(x,x^2,x)=0$} &
\multicolumn{1}{|c|}{$\overline{b}^2\overline{x}\hspace{0.1cm}\overline{a}x\overline{a}xax=
\overline{x}a\overline{b}\hspace{0.1cm}\overline{x}\hspace{0.1cm}\overline{a}xaxb$}
& \multicolumn{1}{|c|}{[1212]} \\
\multicolumn{1}{|c|}{} &
\multicolumn{1}{|c|}{$axbx\overline{b}\hspace{0.1cm}\overline{x}\hspace{0.1cm}\overline{a}b
\overline{x}=xbx\overline{b}x\overline{b}\hspace{0.1cm}\overline{x}\hspace{0.1cm}
\overline{a}^2$} &
\multicolumn{1}{|c|}{[1213]} \\
\multicolumn{1}{|c|}{} &
\multicolumn{1}{|c|}{$\overline{b}a\overline{x}^2bx\overline{a}\hspace{0.1cm}\overline{x}=
\overline{x}\hspace{0.1cm}\overline{b}xa\overline{x}^2b\overline{a}$}
& \multicolumn{1}{|c|}{[1214]} \\ \hline \multicolumn{1}{|c|}{} &
\multicolumn{1}{|c|}{$xa^2x^2bax^2b=axax^2b^2ax^2$} &
\multicolumn{1}{|c|}{[1221]} \\
\multicolumn{1}{|c|}{$(x,x^2,x^2)=0$} &
\multicolumn{1}{|c|}{$\overline{b}^2\hspace{0.1cm}\overline{x}\hspace{0.1cm}\overline{a}x\overline{a}xa
\overline{x}ax=\overline{x}a\overline{b}\hspace{0.1cm}\overline{x}\hspace{0.1cm}\overline{a}
xa\overline{x}axb$} &
\multicolumn{1}{|c|}{[1222]} \\
\multicolumn{1}{|c|}{} &
\multicolumn{1}{|c|}{$axbx\overline{b}\hspace{0.1cm}\overline{x}\hspace{0.1cm}\overline{a}bx\hspace{0.1cm}
\overline{b}\hspace{0.1cm}\overline{x}=xbaxb\hspace{0.1cm}\overline{x}\hspace{0.1cm}
\overline{b}x\overline{b}\hspace{0.1cm}\overline{x}\hspace{0.1cm}\overline{a}$}
& \multicolumn{1}{|c|}{[1223]} \\
\multicolumn{1}{|c|}{} &
\multicolumn{1}{|c|}{$\overline{b}a\overline{x}^2bx\overline{a}\hspace{0.1cm}\overline{b}x^2=
\overline{x}\hspace{0.1cm}\overline{b}a\overline{x}^2ba\hspace{0.1cm}\overline{x}^2b$}
& \multicolumn{1}{|c|}{[1224]} \\ \hline \multicolumn{1}{|c|}{} &
\multicolumn{1}{|c|}{$ax^2bxbx=x^2bax^2b$} &
\multicolumn{1}{|c|}{[2111]} \\
\multicolumn{1}{|c|}{$(x^2,x,x)=0$} &
\multicolumn{1}{|c|}{$xa\overline{x}axb=\overline{x}axbax$} & \multicolumn{1}{|c|}{[2112]} \\
\multicolumn{1}{|c|}{} &
\multicolumn{1}{|c|}{$xb\overline{x}bx\overline{b}\hspace{0.1cm}\overline{x}\hspace{0.1cm}\overline{a}=
axb\overline{x}b\overline{x}b\overline{x}$} &
\multicolumn{1}{|c|}{[2113]} \\
\multicolumn{1}{|c|}{} &
\multicolumn{1}{|c|}{$x\overline{a}\hspace{0.1cm}\overline{b}x^2\overline{a}=
a\overline{x}^2b\overline{a}\hspace{0.1cm}\overline{x}$} &
\multicolumn{1}{|c|}{[2114]} \\ \hline \multicolumn{1}{|c|}{} &
\multicolumn{1}{|c|}{$ax^2bxbax^2=x^2baxax^2b$} &
\multicolumn{1}{|c|}{[2121]} \\
\multicolumn{1}{|c|}{$(x^2,x,x^2)=0$} &
\multicolumn{1}{|c|}{$\overline{x}axba\overline{x}ax=xa\overline{x}a\overline{x}axb$}
& \multicolumn{1}{|c|}{[2122]} \\
\multicolumn{1}{|c|}{} &
\multicolumn{1}{|c|}{$axb\overline{x}b\overline{x}bx=xb\overline{x}baxb\overline{x}$}
& \multicolumn{1}{|c|}{[2123]} \\
\multicolumn{1}{|c|}{} &
\multicolumn{1}{|c|}{$x\overline{a}\hspace{0.1cm}\overline{b}a\overline{x}^2b=
a\overline{x}^2b\overline{a}\hspace{0.1cm}\overline{b}x$} &
\multicolumn{1}{|c|}{[2124]} \\ \hline \multicolumn{1}{|c|}{} &
\multicolumn{1}{|c|}{$ax^2bax^2b^2x=x^2ba^2x^2bxb$} &
\multicolumn{1}{|c|}{[2211]} \\
\multicolumn{1}{|c|}{$(x^2,x^2,x)=0$} &
\multicolumn{1}{|c|}{$\overline{x}\hspace{0.1cm}\overline{a}xa\overline{b}
\hspace{0.1cm}\overline{x}\hspace{0.1cm}\overline{a}xaxb=
\overline{b}\hspace{0.1cm}\overline{x}\hspace{0.1cm}\overline{a}x
\overline{a}\hspace{0.1cm}\overline{x}axbax$} &
\multicolumn{1}{|c|}{[2212]} \\
\multicolumn{1}{|c|}{} &
\multicolumn{1}{|c|}{$xb\overline{x}bx\overline{b}x\overline{b}
\hspace{0.1cm}\overline{x}\hspace{0.1cm}\overline{a}^2=axb\overline{x}
bx\overline{b}\hspace{0.1cm}\overline{x}\hspace{0.1cm}\overline{a}b\overline{x}$}
& \multicolumn{1}{|c|}{[2213]} \\
\multicolumn{1}{|c|}{} &
\multicolumn{1}{|c|}{$x^2\overline{a}\hspace{0.1cm}\overline{b}xa\overline{x}^2b\overline{a}=
a\overline{x}^2ba\overline{x}^2b\overline{a}\hspace{0.1cm}\overline{x}$}
& \multicolumn{1}{|c|}{[2214]} \\ \hline \multicolumn{1}{|c|}{} &
\multicolumn{1}{|c|}{$(ax^2b)^2bax^2=x^2ba(ax^2b)^2$} &
\multicolumn{1}{|c|}{[2221]} \\
\multicolumn{1}{|c|}{$(x^2,x^2,x^2)=0$} &
\multicolumn{1}{|c|}{$\overline{b}\hspace{0.1cm}\overline{x}\hspace{0.1cm}\overline{a}x
\overline{a}\hspace{0.1cm}\overline{x}axba\hspace{0.1cm}\overline{x}ax=
\overline{x}\hspace{0.1cm}\overline{a}xa\hspace{0.1cm}\overline{b}\hspace{0.1cm}\overline{x}
\hspace{0.1cm}\overline{a}xa\hspace{0.1cm}\overline{x}axb$} &
\multicolumn{1}{|c|}{[2222]} \\
\multicolumn{1}{|c|}{} &
\multicolumn{1}{|c|}{$axb\overline{x}bx\overline{b}\hspace{0.1cm}\overline{x}\hspace{0.1cm}\overline{a}bx
\overline{b}\hspace{0.1cm}\overline{x}=xb\overline{x}baxb\overline{x}\overline{b}x\overline{b}
\hspace{0.1cm}\overline{x}\hspace{0.1cm}\overline{a}$} &
\multicolumn{1}{|c|}{[2223]} \\
\multicolumn{1}{|c|}{} &
\multicolumn{1}{|c|}{$(a\overline{x}^2b)^2\overline{a}\overline{b}x^2=
x^2\overline{a}\overline{b}(a\overline{x}^2b)^2$} &
\multicolumn{1}{|c|}{[2224]} \\ \hline
\end{tabular}
\]

\vspace{0.5cm}
\begin{remarks} Let $p,q,r$ be fixed in $\{1,2\}$ and let $n\in\{1,2,3,4\}.$ \begin{enumerate}
\item If $a,b\in\{1,-1\},$ then $\hit_n(a,b)$ is isomorphic to
$\hit_n.$ \item The equation {\bf [pqr3]} follows from {\bf
[rqp2]} by application of the standard involution and replacement
of $(a,b)$ by $(\overline{b},\overline{a}).$ \item The equation
{\bf [2q14]} follows from {\bf [1q24]} by application of the
standard involution and replacement of $(a,b)$ by
$(\overline{b},\overline{a}).\Box$
\end{enumerate}
\end{remarks}

\vspace{0.3cm} \hspace{0.3cm} We shall use the following
equalities: \begin{eqnarray} a &=&
Re(a)+Im(a) \\
Im(a^3) &=& \Big(
3Re(a)^2+Im(a)^2\Big)Im(a)
\end{eqnarray}

\vspace{0.4cm} \hspace{0.3cm} We have the following preliminary
results:

\vspace{0.4cm}
\begin{lemma} $\hit(a,b)$ satisfies to $(x,x^2,x)=0$ if and only if $a,b\in\{1,-1\}.$ In this case $\hit(a,b)$ is
isomorphic to $\hit.$
\end{lemma}

\vspace{0.1cm} {\bf Proof.} The equation {\bf [1211]} gives
$ab=ba,$ by putting $x=a.$ So $Im(a)$ and $Im(b)$ are linearly
dependent. Let now $x$ be a non-zero element in $Im(\hit),$
orthogonal to $a.$ Then, $x$ is also orthogonal to $b,$ and we
have:
\[ xa=\overline{a}x, \hspace{0.2cm} bx=x\overline{b},
\hspace{0.2cm} xa^2=\overline{a}^2x, \hspace{0.2cm}
b^2x=x\overline{b}^2. \] For such an element, the equation {\bf
[1211]} gives \/ $\overline{b}^2=\overline{a}^2.$ As $a$ commutes
with $b,$ we have $a=\pm b$ and equation {\bf [1211]} becomes \/
$[xa^3x,a]=0$ \/ for all $x\in Im(\hit).$ The triple product
identity, applied to $x, a^3,$ gives
\begin{eqnarray} Re(xa^3)[x,Im(a)] &=& 0 \hspace{0.2cm} \mbox{
for all } x\in Im(\hit)
\end{eqnarray}

\vspace{0.2cm} \hspace{0.3cm} If $Im(a)\neq 0,$ then by putting
$x=Im(a)+v,$ where $v$ is a non-zero element in $Im(\hit),$
orthogonal to $a,$ the equalities {\bf (3.4)} and {\bf (3.5)} give
\/ $Im(a^3)=0,$ \/ that is, $a^3=\pm 1.$ Next, we put $x=1+u$
where $u$ is an element in $Im(\hit),$ orthonormal to $a.$ We have
\[ x^2=2u, \hspace{0.2cm} au=u\overline{a}, \hspace{0.2cm}
a^2u=u\overline{a}^2, \hspace{0.2cm} \hspace{0.2cm} aua=u,
\hspace{0.2cm} axa=bxb=a^2+u.
\] The equation {\bf [1211]} gives

\vspace{0.2cm}
\[ (a^2+u)ua^2(1+u)=(1+u)a^2u(a^2+u). \]

It follows that $\overline{a}^2=a^2\in\rit.$ Of this fact
$a=a^3\overline{a}^2\in\rit,$ contradicting the hypothesis
$Im(a)\neq 0.$ So $a=\pm b\in\{1,-1\}$ and $\hit(a,b)$ is
isomorphic to $\hit.\Box$

\vspace{0.2cm}
\begin{lemma} Algebra $^*\hit(a,b)$ satisfies none of four identities $(x,x,x)=0,$
$(x,x^2,x)=0,$ $(x,x,x^2)=0,$ $(x,x^2,x^2)=0.$
\end{lemma}

\vspace{0.1cm} {\bf Proof.} We suppose the opposite.
\begin{enumerate} \item {\bf The first two identities.} We put
$x=a,$ next $x=\overline{b}$ in {\bf [1212]}, and we get
$ab=\overline{b}\hspace{0.1cm} \overline{a}=\overline{b}a.$ It
follows that $a$ and $b$ are scalars, and consequently
$x^2=||x||^2$ for all $x\in\hit.$ This absurdity shows that the
identity $(x,x^2,x)=0,$ and also $(x,x,x)=0,$ are eliminated.
\item {\bf The third identity.} We put $x=1,$ next $x=a$ in {\bf
[1122]}, and we get $\overline{b}a=a^3b=ab.$ Il follows that
$a,b\in\{1,-1\}.$ Thus the equality {\bf [1122]} becomes \/
$||x||^4=||x||^2\overline{x}^2$ for all $x\in\hit.$ Absurd. \item
{\bf The fourth identity.} By putting successively $x=1,$ $x=a,$
$x=a^2$ in equation {\bf [1222]}, we get the equalities
\begin{eqnarray} \overline{b}^2 &=& a\overline{b}ab \\
\overline{b}a &=& ab \\
\overline{b}^2a^2 &=& \overline{a}\overline{b}ab \end{eqnarray}
\end{enumerate}

By taking into account {\bf (3.7)}, the equalities {\bf (3.6)} and
{\bf (3.8)} give $\overline{b}^2=a^2b^2$ and
$\overline{b}^2a^2=b^2,$ respectively. We deduce that
$a^2=b^4\in\{1,-1\}$ and we distinguish following both cases:
\begin{enumerate} \item If $a^2=1,$ the equation {\bf [1222]} becomes $||x||^4\overline{b}^2x=||x||^4\overline{x}$
for all $x\in\hit.$ Absurd. \item If $a^2=-1,$ then by putting
$x=\overline{b}$ in {\bf [1222]} and taking into account {\bf
(3.7)}, we get $b^2=-1.$ Consequently $a^2=b^4=1.$ Absurd.$\Box$
\end{enumerate}

\vspace{0.3cm} \hspace{0.3cm} We get, by means of Remark {\bf 3.5
2)}, a result similar to the Lemma {\bf 3.7}:

\vspace{0.2cm}
\begin{corollary} Algebra $\hit^*(a,b)$ satisfies none of four identities $(x,x,x)=0,$
$(x,x^2,x)=0,$ $(x^2,x,x)=0,$ $(x^2,x^2,x)=0.\Box$
\end{corollary}

\vspace{0.3cm} \hspace{0.3cm} We note $a^{\frac{1}{n}}$ a $n^{th}$
root of $a$ in $\hit.$

\vspace{0.2cm}
\begin{lemma} Algebra $\stackrel{*}{\hit}(a,b)$ satisfies none of four identities $(x,x,x^2)=0,$
$(x,x^2,x^2)=0,$ $(x^2,x,x)=0,$ $(x^2,x^2,x)=0.$
\end{lemma}

\vspace{0.1cm} {\bf Proof.} We suppose the opposite.
\begin{enumerate} \item For the first two identities $(x,x^q,x^2)=0,$ we put $x=1,$
next $x=a^{\frac{1}{2}}$ in {\bf [1q24]}, and we get
$\overline{a}^2=b^2=a.$ We deduct from it that $a$ and $b$ commute
and $a=\pm \overline{b}.$ Thus the equation {\bf [1q24]} becomes
\/ $\overline{b}x^3=\overline{x}\overline{b}^2\overline{x}^2b.$
This gives $b=\pm 1,$ by putting $x=\overline{b}^{\frac{1}{3}},$
and leads to $x^3=\overline{x}^3.$ This absurdity shows that both
equations {\bf [1q24]} cannot come true. \item The Remark {\bf 3.5
3)} shows that the last two equations {\bf [2q14]} cannot come
true also.$\Box$
\end{enumerate}

\vspace{0.3cm} \hspace{0.3cm} We have just proved the following
result:

\vspace{0.3cm}
\begin{proposition} Among those four types of principal isotopes of $\hit$
\begin{enumerate} \item $\hit$ is the unique
algebra satisfying to all identities of the form
$(x^p,x^q,x^r)=0.$ \item Each of three algebras $^*\hit(a,b),$
$\hit^*(a,b),$ $\stackrel{*}{\hit}(a,b),$ of type different from
the first one, can, under certain conditions on the parameters $a$
and $b,$ satisfy both identities \ $(x^2,x,x^2)=(x^2,x^2,x^2)=0.$
Among these three types of algebras
\begin{enumerate} \item $^*\hit$ is the unique
algebra satisfying to $(x^2,x,x)=0,$ or else $(x^2,x^2,x)=0.$
\item $\hit^*$ is the unique satisfying to $(x,x,x^2)=0,$ or else
$(x,x^2,x^2)=0.$ \item $\stackrel{*}{\hit}$ is the unique
satisfying to $(x,x,x)=0,$ or else $(x,x^2,x)=0.\Box$
\end{enumerate}
\end{enumerate}
\end{proposition}

\vspace{0.8cm} \subsection{Specification in case $(x,x^q,x)=0$}

\vspace{0.4cm} \hspace{0.3cm} We are going to establish some
results which will bring supplementary clarification to the
Proposition {\bf 3.9}. First of all a result completing Lemma {\bf
3.6} in the frame of identities of the form $(x,x^q,x)=0:$

\vspace{0.3cm}
\begin{lemma} $\stackrel{*}{\hit}(a,b)$ satisfy to $(x,x^2,x)=0$ if and only if $a,b\in\{1,-1\}.$
In this case $\stackrel{*}{\hit}(a,b)$ is isomorphic to
$\stackrel{*}{\hit}.$
\end{lemma}

\vspace{0.1cm} {\bf Proof.} By putting $x=a$ in {\bf [1214]} we
get $ab=ba.$ Next, we take for $x$ a non-zero element in
$Im(\hit),$ orthogonal to $a.$ The equation {\bf [1214]} give at
first $a^2=b^2,$ that is $b=\pm a$ and consequently $[axa,x]=0$
for all $x\in Im(\hit).$ Lemma {\bf 3.4} and Remark {\bf 3.5 1)}
conclude.$\Box$

\vspace{0.4cm} \hspace{0.3cm} We can state the following result:

\vspace{0.2cm}
\begin{theorem} For a $4$-dimensional absolute-valued algebra $A,$ the following are equivalent:
\begin{enumerate}
\item $A$ is flexible, \item $A$ is third power-associative, \item
$A$ satisfies to the identity $(x,x^2,x)=0,$ \item $A$ is
isomorphic to $\hit$ or $\stackrel{*}{\hit}.$
\end{enumerate}
\end{theorem}

\vspace{0.1cm} {\bf Proof.} The implications {\bf 4 $\Rightarrow$
1 $\Rightarrow$ 2 $\Rightarrow$ 3} are obvious. Besides, the
implication {\bf 3 $\Rightarrow$ 4} is a consequence of Theorem
{\bf 2.7}, the three Lemmas {\bf 3.6, 3.7, 3.11} and Corollary
{\bf 3.8}.$\Box$

\vspace{8cm} \subsection{Specification in cases $(x^p,x^p,x^r)=0$
and $(x^p,x^q,x^q)=0$}

\vspace{0.5cm} \hspace{0.3cm} We have the following result, for
the identity $(x,x,x^2)=0:$

\vspace{0.2cm}
\begin{theorem} For a $4$-dimensional absolute-valued algebra $A,$ the following are equivalent:
\begin{enumerate}
\item $A$ satisfies to the identity $(x,x,x^2)=0,$ \item $A$ is
isomorphic to $\hit$ or $\hit^*.$
\end{enumerate}
\end{theorem}

\vspace{0.1cm} {\bf Proof.} It is enough to show that the first
assertion carry away the second one. Indeed, Lemmas {\bf 3.7} and
{\bf 3.9} show that algebras ${^*\hit}(a,b)$ and
$\stackrel{*}{\hit}(a,b)$ cannot satisfy to $(x,x,x^2)=0.$ For
each one of the two other possibilities $\hit_l(a,b),$ where
$l\in\{1,3\},$ the equation {\bf [112l]} gives $ab=ba,$ by putting
$x=1.$ So $b\in Lin\{1,a\}$ and equations {\bf [1121]}, {\bf
[1123]} give $ax=xa$ for all $x\in Im(\hit),$ that is, $a=\pm 1.$
These equations become $bx^2=x^2b$ and we have $b=\pm 1.$ In this
case $\hit(a,b)$ is isomorphic to $\hit,$ and $\hit^*(a,b)$ is
isomorphic to $\hit^*,$ by virtue of the Remark {\bf 3.5 1)}.
$\Box$

\vspace{0.3cm} \hspace{0.3cm} We deduct from the Theorem {\bf
3.13} and Remark {\bf 3.5 2)} the following result for identity
$(x^2,x,x)=0:$

\vspace{0.2cm}
\begin{theorem} For a $4$-dimensional absolute-valued algebra $A,$ the following are equivalent:
\begin{enumerate}
\item $A$ satisfies to $(x^2,x,x)=0,$ \item $A$ is isomorphic to
$\hit$ or $^*\hit.\Box$
\end{enumerate}
\end{theorem}

\vspace{0.4cm} \hspace{0.3cm} We have the following result, for
identity $(x,x^2,x^2)=0:$

\vspace{0.2cm}
\begin{theorem} For a $4$-dimensional absolute-valued algebra $A,$ the following are equivalent:
\begin{enumerate}
\item $A$ satisfies to $(x,x^2,x^2)=0,$ \item $A$ is isomorphic to
$\hit,$ $\hit^*$ or $\hit^*(1,i).$
\end{enumerate}
\end{theorem}

\vspace{0.1cm} {\bf Proof.} We check easily that algebra
$\hit^*(1,i)$ satisfies to $(x,x^2,x^2)=0$ and it is enough to
show that the first assertion carry away the second one. Indeed,
both Lemmas {\bf 3.7} and {\bf 3.9} show that the algebras
${^*\hit}(a,b)$ and $\stackrel{*}{\hit}(a,b)$ cannot satisfy to
$(x,x^2,x^2)=0.$ For each of two other possibilities
$\hit_l(a,b),$ where $l\in\{1,3\},$ the equation {\bf [122l]}
gives $ab=ba,$ by putting $x=1.$ So $b\in Lin\{1,a\}:=E$ and we
note $E^{\perp}$ the subspace of $\hit$ orthogonal to $E.$ We
distinguish following both cases:
\begin{enumerate} \item If $A=\hit(a,b).$ The same arguments as those in the proof of Theorem {\bf 3.13}
show that $a,b\in\{1,-1\}.$ So $\hit(a,b)$ is isomorphic to
$\hit.$ \item If $A=\hit^*(a,b).$ The equation {\bf [1223]} gives
$b^4=a^2,$ by putting for $x$ a non-zero element in $E^{\perp},$
so $a=\pm b^2.$ The equation {\bf [1223]} becomes
\begin{eqnarray}
b^2xbx\overline{b}\overline{x}\overline{b}x\overline{b}\overline{x}
&=&
xb^3xb\overline{x}\overline{b}x\overline{b}\overline{x}\overline{b}^2
\hspace{0.2cm} \mbox{ for all } x\in\hit
\end{eqnarray} By putting $y=bx,$ the equality {\bf (3.9)} becomes \begin{eqnarray}
by^2\overline{b}\overline{y}\overline{b}y\overline{b}\overline{y}b
&=&
\overline{b}yb^2yb\overline{y}\overline{b}y\overline{b}\overline{y}\overline{b}
\hspace{0.2cm} \mbox{ for all } y\in\hit
\end{eqnarray} We multiply to the right and to the left, by $b,$ the members of the equality {\bf (3.10)} and we
take, for $y,$ a non-zero element in $Im(\hit).$ We get, after a
simplification to the left by $y,$ the equality \
$yby\overline{b}y\overline{b}yb^2=b^2yby\overline{b}y\overline{b}y,$
\ that is
\begin{eqnarray}
[yby\overline{b}y\overline{b}y,b^2] &=& 0 \hspace{0.2cm} \mbox{
for all } y\in Im(\hit)
\end{eqnarray} We compute now, by means of the triple product identity, the expression $yby\overline{b}y\overline{b}y:$
\begin{eqnarray*} yby\overline{b}y\overline{b}y &=& yby\Big( 2(\overline{y}|\overline{b})\overline{b}-\overline{y}\Big)y
\\ &=& yb\Big( 2(\overline{y}|\overline{b})y\overline{b}y-||y||^2y\Big) \\
&=& yb\Bigg( \Big(4(y|b)^2-||y||^2\Big)y-2||y||^2(y|b)b\Bigg) \\
&=& \Big(4(y|b)^2-||y||^2\Big)yby-2||y||^2(y|b)yb^2\\
&=& \Big(4(y|b)^2-||y||^2\Big)\Big( 2(\overline{b}|y)y-||y||^2\overline{b}\Big)-2||y||^2(y|b)yb^2\\
&=& \Big(||y||^2-4(y|b)^2\Big)\Big( 2(b|y)y+||y||^2\overline{b}\Big)-2||y||^2(b|y)yb^2\\
\end{eqnarray*} Now, for $y=Im(b)+u$ where $u\in
E^{\perp}-\{0\},$ we have
\begin{eqnarray*} (y|b) &=& (Im(b)+u|b) \\
&=& (Im(b)|b) \\
&=& ||Im(b)||^2 \end{eqnarray*} and \begin{eqnarray*}
||y||^2-4(y|b)^2 &=& ||Im(b)||^2+||u||^2-4||Im(b)||^4.
\end{eqnarray*} In the other hand, we can choose $u$ rather big in such a way that
\[ ||y||^2-4(y|b)^2:=\lambda\neq 0. \] Thus
\begin{eqnarray*} [yby\overline{b}y\overline{b}y,b^2] &=& \Big[ 2||Im(b)||^2(\lambda
y-||y||^2yb^2)+\lambda||y||^2\overline{b},b^2\Big] \\
&=& 2||Im(b)||^2[\lambda
y-||y||^2yb^2,b^2] \\
&=& 2||Im(b)||^2[\lambda
u-||y||^2ub^2,b^2] \\
&=& 2||Im(b)||^2\Big( \lambda(ub^2-b^2u)-||y||^2(ub^4-b^2ub^2)\Big) \\
&=& 2||Im(b)||^2\Big( \lambda u(b^2-\overline{b}^2)-||y||^2u(b^4-1)\Big) \\
&=& 2||Im(b)||^2u\Big( \lambda(b^2-\overline{b}^2)-||y||^2(b^4-1)\Big) \\
&=& 2||Im(b)||^2u(\lambda\overline{b}^2-||y||^2)(b^4-1).
\end{eqnarray*} Now \begin{eqnarray*} [yby\overline{b}y\overline{b}y,b^2]=0
&\Leftrightarrow&
||Im(b)||^2(\lambda\overline{b}^2-||y||^2)(b^4-1)=0 \\
&\Leftrightarrow& Im(b)=0, \hspace{0.1cm} b^2=\lambda||y||^{-2}
\hspace{0.1cm} \mbox{ or } b^4=1 \\
&\Leftrightarrow& b^4=1.
\end{eqnarray*} If $b=\pm 1,$ then $A=\hit^*(\pm 1,\pm 1)$ is isomorphic to
$\hit^*.$ If $b\in S(Im(\hit)),$ then $b$ is conjugated to $i$ by
Lemma {\bf 2.10} and so $A=\hit^*(\pm 1,b)$ is isomorphic to
$\hit^*(1,i)$ by virtue of Theorem {\bf 2.7}.$\Box$
\end{enumerate}

\vspace{0.3cm} \hspace{0.3cm} We deduce from the Theorem {\bf
3.14} and Remark {\bf 3.5 2)} the following result for the
identity $(x^2,x^2,x)=0:$

\vspace{0.2cm}
\begin{theorem} For a $4$-dimensional absolute-valued algebra $A,$ the following are equivalent:
\begin{enumerate}
\item $A$ satisfy to $(x^2,x^2,x)=0,$ \item $A$ is isomorphic to
$\hit,$ ${^*\hit}$ or ${^*\hit}(i,1).\Box$
\end{enumerate}
\end{theorem}

\vspace{0.8cm} \subsection{Specification in case
$(x^2,x^q,x^2)=0$}

\vspace{0.5cm} \hspace{0.2cm} We are now going to determine those
$4$-dimensional absolute-valued algebras which satisfy to an
identity of the form $(x^2,x^q,x^2)=0.$ We shall see that there is
a variety, among the algebras $\hit_2(a,b)$ and $\hit_3(a,b),$
satisfying such an identity. Let us begin, at first, with the
algebras $\hit_1(a,b)$ and $\hit_4(a,b):$

\vspace{0.2cm}
\begin{lemma} The following are equivalent:
\begin{enumerate}
\item $\hit(a,b)$ satisfy to $(x^2,x,x^2)=0.$ \item $\hit(a,b)$
satisfy to $(x^2,x^2,x^2)=0.$ \item $\hit(a,b)$ is isomorphic to
$\hit.$
\end{enumerate}
\end{lemma}

\vspace{0.1cm} {\bf Proof.} It is enough to show that the third
assertion is consequence of each of both first ones.
\begin{enumerate} \item {\bf 1 $\Rightarrow$ 3.} The equation {\bf
[2121]} gives $abxba=baxab$ for all $x\in\rit\cup Im(\hit).$ So
\begin{eqnarray} abxba &=& baxab \hspace{0.2cm} \mbox{ for all } x\in\hit \end{eqnarray}

\vspace{0.1cm} We deduce from it that
$ba\overline{b}\hspace{0.1cm}\overline{a}=
\overline{b}\hspace{0.1cm}\overline{a}ba$ is a central element of
$\hit,$ and we have \/ $ab=\pm ba.$ Besides, the case $ab=-ba$ is
eliminated at once by making $x=b^{\frac{1}{2}}$ in {\bf [2121]}.
So $ab=ba$ et $b$ belongs in the subalgebra
$\rit_H[a]=Lin\{1,a\}.$ Let now $u$ be an element in $Im(\hit)$
orthonormal to $a,$ then $u$ is orthogonal to $Lin\{1,a\}.$ By
putting $x=1+u$ in {\bf [2121]} and taking into account that the
$\{a,b,\overline{a},\overline{b}\}$ is a commutative set and
$au=u\overline{a}, ub=\overline{b}u,$ we get
$\overline{a}(1+u)b=a(1+u)\overline{b}$ \/ or still
\begin{eqnarray} (1+u)b^2=a^2(1+u)
\end{eqnarray} As $a^2u=u\overline{a}^2,$ we have
$b^2-a^2=u(\overline{a}^2-b^2).$ The trace property of $(.|.),$
and the fact that $b^2-a^2,$ $\overline{a}^2-b^2\in Lin\{1,a\},$
show that $b^2-a^2$ and $u(\overline{a}^2-b^2)$ are orthogonal to
each other. So $b^2=a^2=\overline{a}^2,$ and, taking into account
that $a$ and $b$ commute, we have $b=\pm\overline{a}.$ The Remarks
{\bf 2.8 1), 2.8 3)} and ([CR 08] Corollary {\bf 3.2}) show then
that $\hit(a,b)$ is isomorphic to $\hit.$ \item {\bf 2
$\Rightarrow$ 3.} The equation {\bf [2221]} gives $ab=ba,$ by
taking, for $x,$ a square root of $\overline{a}.$ So $Im(a),$
$Im(b)$ are linearly dependent. By putting $x=1+u$ for norm-one
$u$ in $Im(\hit),$ orthogonal to $a,$ we get \/ $x^2=2u,$
$ua=\overline{a}u$ and $ub=\overline{b}u.$ For such an element
$x,$ the equation {\bf [2221]} gives
$ab=\overline{a}\overline{b}.$ As $ab=ba,$ we have
$b=\pm\overline{a}.\Box$ \end{enumerate}

\vspace{0.2cm}
\begin{lemma} The following are equivalent:
\begin{enumerate} \item $\stackrel{*}{\hit}(a,b)$ satisfy to $(x^2,x,x^2)=0.$
\item $\stackrel{*}{\hit}(a,b)$ satisfy to $(x^2,x^2,x^2)=0.$
\item $\stackrel{*}{\hit}(a,b)$ is isomorphic to
$\stackrel{*}{\hit}.$
\end{enumerate}
\end{lemma}

\vspace{0.1cm} {\bf Proof.} It is enough to show that the third
assertion is consequence of each of both first ones.
\begin{enumerate} \item {\bf 1 $\Rightarrow$ 3.} The equation {\bf
[2124]} gives $x\overline{a}\overline{b}ab=
ab\overline{a}\overline{b}x$ for all $x\in \rit\cup Im(\hit).$ So
\begin{eqnarray} x\overline{a}\overline{b}ab &=&
ab\overline{a}\overline{b}x \hspace{0.2cm} \mbox{ for all }
x\in\hit \end{eqnarray}

It follows easily from ${\bf (3.14)}$ that
$\overline{a}\overline{b}ab=ab\overline{a}\overline{b}\in\rit$ and
then $ab=\pm ba.$ A calculation similar to that of the first part
in the proof of Lemma {\bf 3.17} ends by $b=\pm\overline{a}.$
\item {\bf 2 $\Rightarrow$ 3.} The same as that of the second part
in the proof of Lemma {\bf 3.17}.$\Box$
\end{enumerate}

\vspace{0.5cm} \hspace{0.3cm} We need, for the identity
$(x^2,x,x^2)=0,$ the following supplementary preliminary results:

\vspace{0.2cm}
\begin{lemma} ${^*\hit}(a,b)$ satisfy to $(x^2,x,x^2)=0$ if and only if $a=\pm
1$ and $b\in S(\hit)$ is arbitrary.
\end{lemma}

\vspace{0.1cm} {\bf Proof.} The condition $a=\pm 1$ is sufficient,
according to the equation {\bf [2122]}, and it is enough to show
that it is necessary. Indeed, the equation {\bf [2122]} gives
\begin{eqnarray} [b,(ax)^2] &=& 0 \mbox{ for all } x\in Im(\hit) \end{eqnarray} Next, we put
$x=Im(a)+\overline{a}u,$ where $u$ is an element in $Im(\hit)$
orthogonal to $a,$ and we have \begin{eqnarray} (ax)^2 &=&
Im(a)^4+Re(a)^2Im(a)^2+u^2+2Im(a)^2\Big( Re(a)Im(a)+u\Big)
\end{eqnarray}

\vspace{0.2cm} If $Im(a)\neq 0,$ then the equalities {\bf (3.15)}
and {\bf (3.16)} given $[b,Re(a)Im(a)+u]=0$ for all $u\in
Im(\hit)$ orthogonal to $a.$ Consequently $b$ is a scalar, and the
equation {\bf [2122]} give $[Im(a),x^2]=0$ for all $x\in\hit.$
Absurd.$\Box$

\vspace{0.5cm} \hspace{0.3cm} As a consequence of Lemma {\bf 3.18}
and Remark {\bf 3.5 2)}:

\vspace{0.3cm}
\begin{lemma} $\hit^*(a,b)$ satisfy to $(x^2,x,x^2)=0$ if and only if $b=\pm 1$ and
$a\in S(\hit)$ is arbitrary.$\Box$
\end{lemma}

\vspace{0.3cm} \hspace{0.2cm} We have just proved, by means of
Lemma {\bf 2.10}, the following result:

\vspace{0.2cm}
\begin{theorem} For a $4$-dimensional absolute-valued algebra $A,$ the following are equivalent:
\begin{enumerate}
\item $A$ satisfy to identity $(x^2,x,x^2)=0,$ \item $A$ is
isomorphic to $\hit,$ $\stackrel{*}{\hit},$
${^*\hit}(1,e^{i\alpha})$ or $\hit^*(e^{i\alpha},1)$ with
$\alpha\in [0,\pi[.\Box$
\end{enumerate}
\end{theorem}

\vspace{0.5cm} \hspace{0.3cm} We need, for the identity
$(x^2,x^2,x^2)=0,$ the following supplementary preliminary result:

\vspace{0.3cm}
\begin{lemma} ${^*\hit}(a,b)$ satisfies to $(x^2,x^2,x^2)=0$ if and only if one of the following
conditions holds
\begin{enumerate}
\item $a^2=1$ and $b\in S(\hit).$ In this case ${^*\hit}(a,b)$ is
isomorphic to ${^*\hit}(1,e^{i\alpha})$ where $\alpha\in [0,\pi[.$
\item $a^2=-1$ and $b\in\{1,-1,a,-a\}.$ In this case
${^*\hit}(a,b)$ is isomorphic to ${^*\hit}(i,1)$ or
${^*\hit}(i,i).$
\end{enumerate}
\end{lemma}

\vspace{0.1cm} {\bf Proof.} The equation {\bf [2222]} is obviously
verified if one of both conditions {\bf 1} or {\bf 2} holds. Let
us suppose now that ${^*\hit}(a,b)$ satisfy to $(x^2,x^2,x^2)=0.$
\begin{enumerate} \item We show that {\bf $a^4=1.$} Indeed, the equation {\bf [2222]}
gives at once
\begin{eqnarray} a^2b &=& ba^2
\end{eqnarray} by putting $x=a.$ Now, for norm-one $x$ in $Im(\hit),$
orthogonal to $a$ and $b,$ we have $xa=\overline{a}x$ and
$x\overline{a}=ax.$ The equation {\bf [2222]} gives
$\overline{a}^2b=ba^2$ and the condition $a^4=1$ is obtained from
equality {\bf (3.17)}. \item If $a^2=1,$ then ${^*\hit}(a,b)$ is
isomorphic to ${^*\hit}(1,b)$ by virtue of Remark {\bf 2.8 1)}.
\item If $a^2=-1,$ then the equation {\bf [2222]} is equivalent to
\begin{eqnarray} \overline{b}\overline{x}(axa)\overline{x}(
ax.ba\overline{x}.ax)\in Im(\hit) \hspace{0.2cm} \mbox{ for all }
x\in\hit. \end{eqnarray}
\begin{enumerate} \item Let us show that $ba=\pm ab.$ Indeed, by applying the triple product identity to the
expression $ax(ba\overline{x})ax,$ for an element $x$ in $\hit$
orthogonal to $a,$ we have
\begin{eqnarray*} \overline{x}(axa)\overline{x}\Big(
ax.ba\overline{x}.ax\Big) &=& \overline{x}^3\Big(
2(xa\overline{b}|ax)ax-||x||^2xa\overline{b}\Big)
\\ &=& \overline{x}^3\Big(
2(a\overline{x}\overline{b}|ax)ax-||x||^2a\overline{x}\overline{b}\Big)
\\ &=& \overline{x}^3a\Big(
2(\overline{x}\overline{b}|x)x-||x||^2\overline{x}\overline{b}\Big)
\\ &=& ax^3\Big(x(bx)x\Big) \hspace{0.2cm} \mbox{
by triple product identity }
\\ &=& ax^4bx^2.
\end{eqnarray*}

We note $E$ the subspace of $\hit$ orthogonal to $Lin\{1,a\}.$ By
putting $x=1+u$ where $u\in S(E),$ we get $x^2=2u,$ $x^4=-4.$ The
equation {\bf [2222]} gives $\overline{b}ax^4bx^2\in Im(\hit),$ or
else, $\overline{b}abu\in Im(\hit)$ for all $u\in E.$ The scalar
part $(\overline{b}ab|u)$ of $\overline{b}abu$ vanish for all
$u\in E,$ that is, $\overline{b}ab\in Lin\{1,a\}.$ There are then
scalars $\nu, \rho$ for which $\overline{b}ab=\nu+\rho a,$ or
else, $ab=\nu b+\rho
ba.$ We have \begin{eqnarray*} 0 &=& (a|1)(b|b) \\
&=& (ab|b) \\
&=& (\nu b+\rho ba|b) \\
&=& \nu. \end{eqnarray*} So $ab=\pm ba.$ \item Let us show that
$b\in\{1,-1,a,-a\}.$ Indeed, for any norm-one element $x$ in $E,$
we put
\[ x=1+\frac{a+u}{\sqrt{2}}
\] and we get:
\begin{eqnarray} \overline{x}ax &=& a+u+\sqrt{2}au \end{eqnarray}
and
\begin{eqnarray} \overline{x}ax.a.\overline{x}ax &=&
2a-2u-2\sqrt{2}au
\end{eqnarray}

Now, by putting \[ (y,z)=\Big(
\overline{b}(a-u-\sqrt{2}au)b,-1-au+\sqrt{2}u\Big),
\] and taking into account to the equalities {\bf (3.19)}, {\bf (3.20)} the equation {\bf [2222]} gives
\begin{eqnarray} yz &=& \overline{z}y
\end{eqnarray} If $ba=-ab,$ then $b\in E,$ and we choose $u$ to be orthogonal to $b.$ A first computation gives
\/ $y=-a+u-\sqrt{2}au$ \/ and the equality {\bf (3.21)} leads to
the absurdity $2\sqrt{2}=-2\sqrt{2}.$ So $ab=ba$ and the equality
{\bf (3.21)} gives $\overline{b}^2=b^2,$ or else,
$b\in\{1,-1,a,-a\}.$ Remark {\bf 2.8 1)}, Lemma {\bf 2.10} and
Theorem {\bf 2.7} conclude.$\Box$
\end{enumerate}
\end{enumerate}

\vspace{0.5cm} \hspace{0.3cm} As consequence of Lemmas {\bf 2.10,
3.21} and Remarks {\bf 2.8 1), 3.5 2)}:

\vspace{0.3cm}
\begin{lemma} $\hit^*(a,b)$ satisfy to $(x^2,x^2,x^2)=0$ if and only if one the following conditions holds
\begin{enumerate}
\item $b^2=1$ and $a\in S(\hit).$ In this case $\hit^*(a,b)$ is
isomorphic to $\hit^*(e^{i\alpha},1)$ where $\alpha\in [0,\pi[.$
\item $b^2=-1$ and $a\in\{1,-1,b,-b\}.$ In this case $\hit^*(a,b)$
is isomorphic to $\hit^*(1,i),$ or $\hit^*(i,i).\Box$
\end{enumerate}
\end{lemma}

\vspace{0.3cm} \hspace{0.2cm} We have the following result:

\vspace{0.2cm}
\begin{theorem} For a $4$-dimensional absolute-valued algebra $A,$ the following are equivalent:
\begin{enumerate}
\item $A$ satisfy to $(x^2,x^2,x^2)=0,$ \item $A$ is isomorphic to
\[ \hit, \] \[ {^*\hit}(i,1), {^*\hit}(i,i),
{^*\hit}(1,e^{i\alpha}), \hit^*(1,i), \hit^*(i,i),
\hit^*(e^{i\alpha},1) \mbox{ where } \alpha\in [0,\pi[ \] \[
\mbox{ or } \hspace{0.3cm} \stackrel{*}{\hit}.
\]
\end{enumerate}
\end{theorem}

\vspace{0.1cm} {\bf Proof.} Consequence of Theorem {\bf 2.7} and
Lemmas {\bf 3.16, 3.17, 3.21, 3.22}.$\Box$

\vspace{0.5cm} \hspace{0.3cm} We summarize the results established
in this section as follows:

\vspace{0.2cm}
\begin{theorem} Let $A$ be a $4$-dimensional AAV. Then $A$ satisfy to an identity of the form
$(x^p,x^q,x^r)=0,$ for fixed integers $p,q,r\in\{1,2\},$ if and
only if $A$ is isomorphic to one of principal isotopes of $\hit$
described according to the following table:

\vspace{0.2cm}
\[ \begin{tabular}{cc} \\ \hline
\multicolumn{1}{|c|}{$A$ satisfy to} & \multicolumn{1}{|c|}{The list of isomorphism classes} \\
\hline \multicolumn{1}{|c|}{$(x,x,x)=0$ or $(x,x^2,x)=0$} &
\multicolumn{1}{|c|}{$\hit, \stackrel{*}{\hit}$} \\
\hline \multicolumn{1}{|c|}{$(x,x,x^2)=0$} &
\multicolumn{1}{|c|}{$\hit,$ $\hit^*$} \\
\hline \multicolumn{1}{|c|}{$(x,x^2,x^2)=0$} &
\multicolumn{1}{|c|}{$\hit,$ $\hit^*,$ $\hit^*(1,i)$} \\
\hline \multicolumn{1}{|c|}{$(x^2,x,x)=0$} &
\multicolumn{1}{|c|}{$\hit,$ ${^*\hit}$}
\\ \hline \multicolumn{1}{|c|}{$(x^2,x^2,x)=0$} &
\multicolumn{1}{|c|}{$\hit,$ ${^*\hit},$ ${^*\hit}(i,1)$}
\\ \hline \multicolumn{1}{|c|}{$(x^2,x,x^2)=0$} &
\multicolumn{1}{|c|}{$\hit,$ $\stackrel{*}{\hit}$} \\
\multicolumn{1}{|c|}{} &
\multicolumn{1}{|c|}{${^*\hit}(1,e^{i\alpha}),$
$\hit^*(e^{i\alpha},1)$ with $\alpha\in [0,\pi[$}
\\ \hline \multicolumn{1}{|c|}{} &
\multicolumn{1}{|c|}{$\hit, \stackrel{*}{\hit}$}
\\ \multicolumn{1}{|c|}{$(x^2,x^2,x^2)=0$} &
\multicolumn{1}{|c|}{${^*\hit}(i,1),$ ${^*\hit}(i,i),$
$\hit^*(1,i),$
$\hit^*(i,i)$} \\
\multicolumn{1}{|c|}{} &
\multicolumn{1}{|c|}{${^*\hit}(1,e^{i\alpha}),$
$\hit^*(e^{i\alpha},1)$ with $\alpha\in [0,\pi[$}
\\ \hline
\end{tabular}
\]
\end{theorem}

\vspace{0.2cm}
\begin{remarks}. \begin{enumerate}
\item All isomorphism classes, in above table, are represented by
commuting parameters $a,b.$ Theorem {\bf 2.9} shows that every
$4$-dimensional absolute-valued algebra satisfying to an identity
of the form $(x^p,x^q,x^r)=0,$ for fixed $p,q,r$ in $\{1,2\},$
contains $2$-dimensional subalgebras. Concretely, algebras $\hit$
and $\stackrel{*}{\hit}$ contain $\cit$ and $\stackrel{*}{\cit},$
respectively. Those algebras ${^*\hit}(a,b)$ and $\hit^*(a,b),$ in
above table, contain $^*\cit$ and $\cit^*,$ respectively. \item
Among those $4$-dimensional absolute-valued algebras, each
identity of the form $(x^p,x^q,x^r)=0,$ other than
$(x^2,x^2,x^2)=0,$ carry away $(x^2,x^2,x^2)=0.$ The converse is
false. \item Every $4$-dimensional absolute-valued algebra
satisfying any one of identities \ $(x,x,x)=0,$ $(x,x^2,x)=0,$
$(x,x,x^2)=0$ or $(x^2,x,x)=0$ \ cannot be single-generated.$\Box$
\end{enumerate}
\end{remarks}

\vspace{1cm} \section{Variety in dimension eight for
absolute-valued algebras satisfying $(x^2,y,x^2)=0$}

\vspace{0.8cm} \subsection{Possibility of absence of
$4$-dimensional subalgebras}

\vspace{0.6cm} \hspace{0.3cm} In this section we are going to put
evidence the unlimitedness of the class of $8$-dimensional
absolute-valued algebras satisfying to $(x^2,y,x^2)=0$ by giving
examples of such an algebras which contain no $4$-dimensional
subalgebras.

\vspace{0.3cm} \hspace{0.3cm} Gleichgewicht [G 63] defined on the
normed space, underlying to an absolute-valued algebra $(A,||.||)$
with involution $*,$ a new product $\odot$ by putting \ $x\odot
y=x^*y$ for all $x,y\in A.$ This provides a new absolute-valued
algebra $(A,\odot,||.||):=A_*$ called a cracovian algebra
generated by $A.$ Besides, there exists an element $e\in A,$ such
that $x^*x=||x||^2e$ for all $x\in A$ [G 63]. This shows that $e$
is the unique non-zero idempotent for cracovian algebra $A_*.$

\vspace{0.3cm} \hspace{0.3cm} Now, let us remind that every
absolute-valued algebra $(A,||.||)$ containing a non-zero central
idempotent is endowed with an involution [Elm 90] and that the
norm $||.||$ comes from an inner product $(.|.)$ [U 61]. It is
shown that if $(A,||.||,(.|.))$ is an absolute-valued real algebra
which is not isomorphic to $\stackrel{*}{\cit},$ containing a
non-zero central $e,$ admits one and only one involution $*$ ([RR
09] Corollary {\bf 3.8}), given by $x\mapsto x^*:=2(x|e)e-x$ ([RR
09] Proposition {\bf 3.5}). It follows immediately the result:

\vspace{0.2cm} \begin{proposition} Let $A$ be an absolute valued
real algebra which is not isomorphic to $\stackrel{*}{\cit},$
containing a non-zero central idempotent $e,$ and let $*$ its
unique involution. Then every subspace of $A$ containing $e$ is
$*$-invariant.$\Box$
\end{proposition}

\vspace{0.3cm} \hspace{0.3cm} For every linear isometry $f$ of
euclidian space $\oit,$ fixing $1,$ the normed space $\oit$
endowed with the product $\bullet$ defined, for all $x,y\in\oit,$
by $x\bullet y=f(x)f(y)$ becomes an absolute-valued real algebra
$\oit^f$ with central idempotent $1$ [RR 09]. According to ([RR
09] Proposition {\bf 3.5}) the mapping \/
$\alpha.1+u\stackrel{*}{\mapsto}\alpha.1-u$ \/ is an involution on
algebra $\oit^f=\rit.1\oplus Im(\oit).$ In the other hand, it is
well known that any finite-dimensional absolute-valued algebra
contains a non-zero idempotent [Se 54]. With the same notations as
above, we have the following result:

\vspace{3cm}
\begin{proposition} Every subalgebra $B$ of $(\oit^f)_*$ is a subalgebra of $\oit^f.$
\end{proposition}

\vspace{0.1cm} {\bf Proof.} The unique non-zero idempotent $e=1$
of algebra $(\oit^f)_*$ belongs to $B.$ So the subspace $B$ of
$\oit$ is $*$-invariant by virtue of Proposition {\bf 4.1}, and
then $B$ is a subalgebra of $\oit^f.\Box$

\vspace{0.5cm} \hspace{0.3cm} We have the following result about a
flexible idempotent:

\vspace{0.2cm} \begin{proposition} Let $A$ be an absolute-valued
algebra $A$ of finite-dimension $n\geq 2,$ having a nonzero
flexible idempotent $e.$ Then $A$ contains a $2$-dimensional
subalgebra invariant under $L_e$ and $R_e.$
\end{proposition}

\vspace{0.1cm} {\bf Proof.} The operators $L_e,R_e$ are linear
isometries fixing $e,$ and induce isometries on the orthogonal
space $(\rit e)^{\perp}:=E.$ As $E$ has odd dimension and $L_e$
commutes with $R_e,$ there exists common norm-one eigenvector
$u\in E$ for both $L_e$ and $R_e$ associated to eigenvalues
$\alpha, \beta\in\{1,-1\}.$ Now, the normed space $A$ endowed with
multiplication \ $x\odot y=L_e(x)R_e(y)$ \ is an absolute-valued
algebra for which we check easily that $e$ is a central
idempotent. So \ $u\odot u=-e$ ([Elm 90] Lemma {\bf 3.3}), that
is, \ $u^2=-\alpha\beta e.$ It follows that $Lin\{e,u\}$ is a
$2$-dimensional subalgebra of $A$ invariant under both $L_e$ and
$R_e.\Box$

\vspace{0.4cm}
\begin{remark} There exists an appropriate linear isometry $f,$ of euclidian space $\oit,$ for which the algebra
$\oit^f$ contains no $4$-dimensional subalgebras {\em ([RR 09]
Theorem {\bf 4.8}, Corollary {\bf 4.9})}.$\Box$
\end{remark}

\vspace{0.4cm} \hspace{0.3cm} A finite-dimensional algebra $A$
over a field $K$ is said to be of degree $n\in\nit$ if $n$ is the
minimum natural number such that all single-generated subalgebras
of $A$ have dimension $\leq n.$ It follows from Theorem {\bf 2.5},
that finite-dimensional absolute valued real algebras have degree
$1,2,4$ or $8.$ Rodr\'iguez ([Rod 94] Theorem {\bf 3}) showed that
\ $\cit,$ ${^*}\cit,$ $\cit^*,$ $\stackrel{*}{\cit},$ $\hit,$
${^*}\hit,$ $\hit^*,$ $\stackrel{*}{\hit},$ $\oit,$ ${^*}\oit,$
$\oit^*,$ $\stackrel{*}{\oit}$ (the standard isotopes of $\cit,$
$\hit,$ $\oit$) and $\pit$ (the pseudo-octonion algebra [Ok 78])
are the only absolute valued real algebras of degree $2.$

\vspace{4cm}
\begin{remarks} . \begin{enumerate}
\item Every finite-dimensional third power-associative absolute
valued real algebra has degree $\leq 2$ {\em ([Elm 87], [EP 94])}.
\item Every $8$-dimensional absolute valued real algebra of degree
$2$ contains $4$-dimensional subalgebras.$\Box$
\end{enumerate}
\end{remarks}

\vspace{0.4cm} \hspace{0.3cm} Now, with the same notation as in
Remark {\bf 4.4}, we state the following main result of this
subsection:

\vspace{0.5cm}
\begin{proposition} Let $B$ denotes algebra $(\oit^f)_*.$ Then
\begin{enumerate} \item $B$ contains a $2$-dimensional subalgebra but no $4$-dimensional
subalgebras. As immediate consequence, $B$ has degree $8.$ \item
$B$ satisfy to \/ $(x^2,y,x^2)=0$ but do not satisfy to \
$(x,x^2,x)=0.$
\end{enumerate}
\end{proposition}

\vspace{0.1cm} {\bf Proof.} \begin{enumerate} \item A direct
computation shows that $e=1\in B$ is a flexible idempotent, so $B$
contains a $2$-dimensional subalgebra according to Proposition
{\bf 4.3}. In other hand $B$ contains no $4$-dimensional
subalgebras by virtue of Remark {\bf 4.4} and Proposition {\bf
4.2}. By Remark {\bf 4.5.2} $B$ has degree $8.$ \item Follows by
direct computation.$\Box$
\end{enumerate}

\vspace{0.8cm} \subsection{On $8$-dimensional absolute-valued
algebras satisfying $(x,x^2,x)=0$}

\vspace{0.6cm} \hspace{0.3cm} Among those $8$-dimensional
absolute-valued algebras satisfying $(x,x^2,x)=0$ we give
sufficient conditions which insure third power-associativity:

\vspace{4cm} \begin{theorem} Let $A$ be a finite-dimensional
absolute-valued algebra satisfying $(x,x^2,x)=0,$ then the degree
of $A$ is $\neq 4.$ Moreover, the following are equivalent:
\begin{enumerate} \item $A$ is flexible, \item $A$ is third power-associative, \item $A$
has degree $\leq 2,$ \item $A$ has degree $\neq 8,$ \item $A$ is
isomorphic to \ $\rit,$ $\cit,$ $\stackrel{*}{\cit},$ $\hit,$
$\stackrel{*}{\hit},$ $\oit,$ $\stackrel{*}{\oit},$ or $\pit.$
\end{enumerate}
\end{theorem}

\vspace{0.1cm} {\bf Proof.} The first proposition is a consequence
of Remark {\bf 3.25.3}. This shows the equivalence {\bf 3
$\Leftrightarrow$ 4}. In the other hand, the fifth assertion carry
away all four above ones. Moreover {\bf 1 $\Rightarrow$ 2} is
clear. Now, the implication {\bf 2 $\Rightarrow$ 3} follows from
Remark {\bf 4.5.1}. Finally {\bf 3 $\Rightarrow$ 5} follows by
computation using [Rod 94]. This completes the proof.$\Box$

\vspace{0.3cm}
\begin{problem} Let $A$ be an $8$-dimensional absolute-valued real algebra satisfying both
identities \ $(x,x^2,x)=(x^2,y,x^2)=0.$ It is interesting to know
if algebra $A$ is third power associative, or even contains
$4$-dimensional subalgebras.$\Box$
\end{problem}

\vspace{0.2cm} \hspace{0.3cm} {\scriptsize\bf D\'EPARTEMENT DE
MATH\'EMATIQUES ET INFORMATIQUE, \\ FACULT\'E DES SCIENCES BEN
M'SIK, UNIVERSIT\'E HASSAN II-MOHAMMEDIA, \\ B.P. 7955 CASABLANCA
(MOROCCO)}

\vspace{0.1cm} \hspace{0.3cm} {\scriptsize\bf E-mail address:
ahchandid@hotmail.com}

\vspace{0.3cm} \hspace{0.3cm} {\scriptsize\bf DEPARTAMENTO DE
\'ALGEBRA Y AN\'ALISIS MATEM\'ATICO, \\ UNIVERSIDAD DE ALMER\'IA.
04071 ALMER\'IA (SPAIN)}

\vspace{0.1cm} \hspace{0.3cm} {\scriptsize\bf E-mail address:
mramirez@ual.es}

\vspace{0.3cm} \hspace{0.3cm} {\scriptsize\bf D\'EPARTEMENT DE
MATH\'EMATIQUES ET INFORMATIQUE, \\ FACULT\'E DES SCIENCES BEN
M'SIK, UNIVERSIT\'E HASSAN II-MOHAMMEDIA, \\ B.P. 7955 CASABLANCA
(MOROCCO)}

\vspace{0.1cm} \hspace{0.3cm} {\scriptsize\bf E-mail address:
abdellatifro@yahoo.fr}

\end{document}